\documentclass[a4paper,12pt]{amsart}

\usepackage[english]{babel}
\usepackage{amsthm}
\usepackage{amssymb}

\newcommand{\Set}[1]{\left\{\, #1 \,\right\}}
\newcommand{\Span}[1]{\langle\, #1 \,\rangle}
\newcommand{\Order}[1]{\lvert #1 \rvert}

\DeclareMathOperator{\Sym}{Sym}
\DeclareMathOperator{\Alt}{Alt}
\DeclareMathOperator{\AGL}{AGL}
\DeclareMathOperator{\PSL}{PSL}
\DeclareMathOperator{\GL}{GL}
\DeclareMathOperator{\SL}{SL}
\DeclareMathOperator{\GF}{GF}

\DeclareMathOperator{\Soc}{Soc}
\DeclareMathOperator{\Aut}{Aut}
\DeclareMathOperator{\Out}{Out}

\renewcommand{\phi}[0]{\varphi}
\renewcommand{\theta}[0]{\vartheta}
\renewcommand{\epsilon}[0]{\varepsilon}

\newtheorem{dummy}{Dummy}
\numberwithin{equation}{section}

\newtheorem{theorem}[dummy]{Theorem}

\newtheorem*{Assu}{Cryptographic Assumptions}
\newcommand{\theAssu}{Cryptographic Assumptions}

\begin{document}

\bibliographystyle{amsalpha}

\date{15 March 2009
}

\title[O'Nan-Scott and cipher groups]
{An application of the O'Nan-Scott theorem to the group generated by
the round functions of an AES-like cipher} 

\author{A.~Caranti}

\address[A.~Caranti]{Dipartimento di Matematica\\
  Universit\`a degli Studi di Trento\\
  via Sommarive 14\\
  I-38050 Povo (Trento)\\
  Italy} 

\email{andrea.caranti@unitn.it} 

\urladdr{http://science.unitn.it/\~{ }caranti/}

\author{F.~Dalla Volta}

\address[F.~Dalla Volta]{Dipartimento di Matematica e Applicazioni\\
  Edificio U5\\
  Universit\`a degli Studi di Milano--Bicocca\\
  via Roberto Cozzi 53\\
  I-20125 Milano\\
  Italy}

\email{francesca.dallavolta@unimib.it}

\urladdr{http://www.matapp.unimib.it/\~{ }dallavolta/}

\author{M.~Sala}

\address[M.~Sala]{Dipartimento di Matematica\\
  Universit\`a degli Studi di Trento\\
  via Sommarive 14\\
  I-38050 Povo (Trento)\\
  Italy} 

\email{sala@science.unitn.it}

\begin{abstract}
  In  a previous  paper,  we  had proved  that  the permutation  group
generated   by  the  round   functions  of   an  AES-like   cipher  is
primitive. Here  we apply the O'Nan Scott  classification of primitive
groups to  prove that this group  is the alternating 
group.
\end{abstract}


\thanks{First author partially supported by MIUR-Italy via PRIN 2006014340-002
 ``Lie  algebras  and rings.  Groups.  Cryptography''.  Second  author
  partially supported  by MIUR-Italy via PRIN 2007  ``Group theory and
  applications''.}

\maketitle

\thispagestyle{empty}

\section{Introduction}

According  to  Shannon~\cite[p.~657]{Sha49},  a  cipher  ``is  defined
abstractly   as   a  set   of   transformations''.   Coppersmith   and
Grossman~\cite{CopGro},  and later  in  1988 Kaliski,  Rivest and 
Sherman~\cite{DES},  called  attention to  the  group  generated by  a
cipher.  One of the motivations for the work of Kaliski et al.\ is
that at that time Triple DES  was being suggested as an improvement to
DES. This  meant replacing  the use of  single DES   transformation
$T_{a}$, where $a$ is a key, with the composition $T_{a} T_{b} T_{c}$,
where $a, b,   c$ are three DES  keys. If it was the  case that the
transformations of  DES form a group,  then Triple   DES would have
been of  course no  more than DES  itself. More generally,  Kaliski et
al.\ showed that  if  the group generated by the transformations of
a  cipher is  too small,  then the  cipher is  exposed   to certain
cryptanalytic attacks.

It was later proved by Wernsdorf~\cite{DESW1} that the group generated
by the  round functions  of DES (which  are even permutations)  is the
alternating group.  This implies that  the group generated by  the DES
transformations  with  independent  subkeys  is also  the  alternating
group. (We  are not aware  of any work  in this context that  tries to
take account of the key schedule.)

Wernsdorf  used  ad  hoc  methods  in~\cite{WAES} to  prove  that  the
permutation group $G$  generated by the round functions  of AES is the
alternating group. (Here, too, these functions are even permutations.)
Sparr  and Wernsdorf  have recently  given another,  permutation group
theoretic proof in~\cite{SpWe}.

The goal  of this  paper is to  give a  different proof of  this fact,
building upon our earlier  paper~\cite{Impr}. There we had proved that
the group $G$ is primitive. In  the course of doing that we answered a
question of  Paterson~\cite{Pat} about the possibility  of embedding a
trapdoor in a  cipher by having the group generated  by the cipher act
imprimitively.

In  this paper we  work under  certain cryptographic  assumptions (see
Section~\ref{sec:prelim}) that are a stripped down, simplified version
of those of~\cite{Impr}. (These are  also satisfied by AES.)  We first
give,  for the  convenience  of the  reader,  a short  group-theoretic
version of the main  result of~\cite{Impr} under these assumptions. We
then appeal  to the O'Nan-Scott classification of  primitive groups to
prove  that  the   group  generated  by  the  round   functions  of  a
cryptosystem satisfying our assumptions is the alternating group.

We are very grateful to Ralph Wernsdorf for several useful suggestions.

\section{Preliminaries}\label{sec:prelim}

In the rest of the paper, we tend to adopt the notation of~\cite{AES}.

Let $V  = V(d, 2)$, the vector  space of dimension $d$  over the field
$\GF(2)$ with two  elements, be the state (or  message) space. $V$ has
$n = 2^{d}$ elements.

For any $v \in V$, consider the translation by $v$, that is the map
\begin{equation*}
  \begin{aligned}
  \sigma_{v} : V &\to V,\\
  w &\mapsto w + v.
  \end{aligned}
\end{equation*} 
In particular, $\sigma_{0}$ is the identity map on $V$.  The set
\begin{equation*}
  T = \Set{\sigma_{v} : v \in V}
\end{equation*}
is an elementary abelian, regular subgroup of $\Sym(V)$. In fact, the map
\begin{equation}\label{eq:isoVT}
  \begin{aligned}
    V &\to T\\ v &\mapsto \sigma_{v}\\
  \end{aligned}
\end{equation}
is an isomorphism of the additive group $V$ onto the multiplicative
group $T$.

We consider  a \emph{key-alternating block  cipher} (see Section~2.4.2
of~\cite{AES}) which  consists of  a fixed number  of iterations  of a
function  of the  form  $\rho \sigma_{k}$,  where  $k \in  V$. Such  a
function is called  a \emph{round function}, and the  parameter $k$ is
called the \emph{round key}.   (We write maps left-to-right, so $\rho$
operates first.)  Here $\rho$ is  a fixed permutation operating on the
vector space $V$.  Therefore each  round consists of an application of
$\rho$, followed by a key addition.  This covers for instance AES with
\emph{independent} subkeys.   Let $G =  \Span{\rho \sigma_{k} :  k \in
  V}$  be the  group of  permutations of  $V$ generated  by  the round
functions. Choosing $k = 0$ we see  that $\rho \in G$, and thus $T \le
G$. It follows that $G = \Span{T, \rho}$.

We assume $\rho = \gamma
\lambda$, where $\gamma$ and $\lambda$ are permutations. Here $\gamma$
is a bricklayer transformation, consisting of a number of S-boxes. The
message space $V$ is written as a direct sum
\begin{equation*}
  V = V_{1} \oplus \dots \oplus V_{n_{t}},
\end{equation*}
$n_{t} > 1$, where each $V_{i}$ has the same dimension $m>1$ over $\GF(2)$. As $n_{t} > 1$, this implies that $d=m n_t$ is not a prime number. For $v
\in V$, we will write $v = v_{1} + \dots + v_{n_{t}}$, where $v_{i} \in
V_{i}$. Also, we consider the projections $\pi_{i} : V \to V_{i}$,
which map $v \mapsto v_{i}$. We have
\begin{equation*}
  v \gamma = v_{1} \gamma_{1} \oplus \dots \oplus v_{n_{t}} \gamma_{n_{t}},
\end{equation*}
where the $\gamma_{i}$ are S-boxes, which we allow to be different for
each $V_{i}$.

$\lambda$ is a linear function (usually called a linear mixing
layer). The only assumption we will be making about $\lambda$ is
Cryptographic Assumption~\eqref{item:affine} below.

In AES  the S-boxes  are all  equal, and consist  of inversion  in the
field  $\GF(2^{8})$   with  $2^{8}$   elements  (see  later   in  this
paragraph), followed  by an affine  transformation, that is,  a linear
transformation, followed  by a  translation. When interpreting  AES in
our scheme, we take advantage  of the well-known possibility of moving
the  linear part  of the  affine transformation  to the  linear mixing
layer, and incorporating the translation  in the key addition (see for
instance~\cite{BES}). Thus in  our scheme for AES we have  $m = 8$, we
identify each $V_{i}$  with $\GF(2^{8})$, and we take  $x \gamma_{i} =
x^{2^{8}-2}$,  so that  $\gamma_{i}$  maps nonzero  elements to  their
inverses,  and  zero  to  zero.  As  usual we  will  simply  say  that
$\gamma_{i}$ acts by inversion.

We will work under the following
\begin{Assu}
Consider an AES-like cryptosystem  as described above, which satisfies
the following conditions.
\begin{enumerate}
\item  
  \label{item:01} 
  $0 \gamma = 0$ and $\gamma^{2} = 1$, the identity transformation.
\item There is $1 \le r < m/2$ such that the following hold.
  \begin{enumerate}
  \item \label{item:r} 
    For all $0 \ne v \in V_{i}$, the image of the map $V_{i}
    \to V_{i}$,  which maps  $x \mapsto (x  + v) \gamma_{i}  + x
    \gamma_{i}$, has size greater than $2^{m-r-1}$, and it is not a
    coset of a subspace.
  \item \label{item:b}There is no subspace of $V_{i}$, invariant
    under $\gamma_{i}$, of codimension less than or equal to $2 r$.
  \end{enumerate}
  \item 
    \label{item:affine}
    There are no subspaces $U,  U', U''$ (except $\Set{0}$ and $V$) that
    are the sum of some of the $V_{i}$, and such that $U \lambda = U'$
    and $U' \lambda = U''$.
  \end{enumerate}
\end{Assu}

In~\cite{Impr}  we  have proved  under  certain  abstract and  general
assumptions a result that specializes to the following:
\begin{theorem}\label{th:primitive}
Suppose  a cryptosystem  satisfies  the~\theAssu. Then  the group  $G$
generated by its round functions acts primitively on the message space
$V$.
\end{theorem}
We    give   a    short,    group-theoretic   proof    of   this    in
Section~\ref{sec:referee}.   This we  do  for the  convenience of  the
reader,  as  we   will  need  to  refer  to  part   of  the  proof  in
Section~\ref{sec:wreath}. We  are grateful  to the referee  of another
paper for this proof.

In the rest of the paper we prove the following
\begin{theorem}\label{uth}
Suppose a cryptosystem satisfies the~\theAssu. 

Then  the  group  $G$  generated  by its  round  functions  is  the alternating group
$\Alt(V)$.

The  same holds  for  the  group generated  by  the cryptosystem  with
independent subkeys.
\end{theorem}

A word about the parity of the  group $G$ is in order here.  Over $V =
V(d,2)$,  non-trivial  translations  are clearly  involutions  without
fixed points, and thus even permutations.  Also, for $d > 2$ the group
$\GL(d,2)  = \SL(d,2)$ is  perfect, so  that in  particular it  has no
(normal) subgroup of order $2$, and it is thus contained in $\Alt(V)$.

We now  show that $\gamma$ is also  even, so that $G  \le \Alt(V)$. In
fact, $\gamma$ is  the product of $n_t$ permutations  $g_i$, acting as
$\gamma _i$  on $V_i$, and as  the identity on $V_j$,  $j\neq i$. This
means that every $2$-cycle  in $\gamma_i$ gives rise to $2^{d-m}$
$2$-cycles in $g_{i}$. Now the number  $2^{d-m}$ is even, as $d - m
= n_{t}  m - m  > m$,   $n_{t} >  1$ by assumption,  and $m >  2$ by
Cryptographic Assumption~\eqref{item:01}. It follows that each $g_{i}$
is even, and thus so is $\gamma$.  (The same argument proves that $\gamma $
is even,  even without  assuming that  it is an  involution, as  we do
here.)

Condition~\eqref{item:01}  is clearly  satisfied  by AES.  As we  said
above,  we take  advantage here  of the  possibility of  assuming that
$\gamma$ is simply componentwise inversion.

Condition~\eqref{item:r} is also well-known to be
satisfied,  with  $r =  1$  (see~\cite{Ny}  but also~\cite{DR}), as
the image of that map has size $2^{7}-1$.

As to Condition~\eqref{item:b}, it is  also satisfied by AES with $r =
1$. For  that, one could  just use \textsf{GAP}~\cite{GAP4}  to verify
that the  only nonzero subspaces  of $\GF(2^{8})$ which  are invariant
under inversion are the subfields.   However, this can also be derived
from  a  more general  result  of~\cite{GGSZ}  and \cite{Matt},  which
states that the only nonzero additive subgroups of $\GF(2^{m})$, which
contain  the  inverse  of  all  of their  nonzero  elements,  are  the
subfields.

Condition~\eqref{item:affine} follows from  the   properties  of  the
components       \texttt{MixColumns}       \cite[3.4.3]{AES}       and
\texttt{ShiftRows} \cite[3.4.2]{AES} of the linear mixing layer (which
are not altered by the fact that we have incorporated in it the linear
part of the  S-boxes).  In fact, suppose, without  loss of generality,
that $U \supseteq  V_{1}$.  Then $U'$ contains the  whole first column
of the state, and $U'' = V$, a contradiction. This argument is a
vestigial form of the Four-Round Propagation Theorem~\cite[9.5.1]{AES}.

\section{Primitivity}\label{sec:referee}

In this section we give a proof of Theorem~\ref{th:primitive}.

Suppose  for   a  contradiction  that   $G=\langle  T,\rho\rangle$  is
imprimitive on $V$,  so that any block system for $G$  is given by the
cosets of some  subspace $U$ of $V$. This is because,  as it is proved
in~\cite{Impr}, a block system for $G$  is also a block system for the
group $T$ of translations.

Now $\rho  = \gamma \lambda$, with  $\lambda$ linear, and  $0 \gamma =
0$. Thus  $U \rho  = U$, and  $U' =  U \gamma =  U \lambda^{-1}$  is a
subspace.

Suppose  firstly  that  $U=V_{i_1}\oplus\cdots  \oplus V_{i_l}$  is  a
direct  sum of  some of  the  subspaces $V_i$  ($l<n_t$). Then,  $U'=U
\gamma =U$, so that $U'=U$ is $\lambda$-invariant; this contradicts
Cryptographic Assumption~\eqref{item:affine}.

Thus there exists $i$ such  that $U \not\supseteq V_{i}$, but there is
$u  \in  U$, such  that  its $i$-th  component  $u_{i}  \in V_{i}$  is
nonzero. We  claim that  $U \cap  V_{i}$ is nonzero.  Take any  $v \in
V_{i}$. Then $(u +  v) \gamma + v \gamma \in U'$,  so that $u \gamma +
(u +  v) \gamma + v  \gamma \in U'$.  The latter element has  all zero
components, expect possibly the $i$-th one, which is $u_{i} \gamma_{i}
+ (u_{i} + v) \gamma_{i} + v  \in U' \cap V_{i}$. Were the latter zero
for all  $v \in V_{i}$,  then the map  $V_{i} \to V_{i}$ that  maps $v
\mapsto (u_{i} + v) \gamma_{i} + v \gamma_{i}$ would be constant, thus
contradicting Cryptographic Assumption~\eqref{item:r}.

Thus   there  exists   $i$  such   that  both   $U_i=U\cap   V_i$  and
$U'_i=(U_i)\gamma  _i=U'\cap  V_i$ are  nonzero,  proper subspaces  of
$V_i$ of the same dimension, and
$$\gamma _i:  {V_i}/{U_i}\rightarrow {V_i}/{U'_i}.$$ If  $x \in
V_{i}$, and $v\in U_i, \,  
v\neq  0$, then  $x+v$ and  $x$ are  in the  same coset  of  $U_i$, so
$(x+v)\gamma  _i$  and   $x\gamma  _i$  are  in  the   same  coset  of
$U'_i$. Thus the set
\begin{equation*}
\{ (x+v)\gamma _i+x\gamma _i : x\in V_i\}
\end{equation*}
is a subset of  $U'_i$, and by Cryptographic Assumption~\eqref{item:r}
$U_i$ and  $U'_i$ have size greater  than $2^{m-r-1}$, that  is to say
dimension at least $m-r$ or  equivalently codimension at most $r$. The
codimension of $U_i \cap U'_i$ is therefore at most $2r$, so $U_i \cap
U'_i$  cannot  be  $\gamma  _{i}$-invariant because  of  Cryptographic
Assumption~\eqref{item:b}. This means there exists $z\in U_i\cap U'_i$
such that $z\gamma _i\notin  U_i\cap U'_i$, so $z\gamma _i\notin U_i$,
as $z\gamma _i\in  U'_i$. However, $U'_i$ is the  image of $U_i$ under
the  bijective map  $\gamma _i$,  so $z=z\gamma  _i^2\notin  U'_i$, as
$z\gamma  _i\notin  U_i$. Thus  $z\notin  U_i\cap  U'_i$,  which is  a
contradiction.

\section{O'Nan-Scott}

In  this  section  we  prove  Theorem~\ref{uth}. We  first  state  the
O'Nan-Scott  classification of  primitive groups  for the  case  of the
maximal primitive subgroups of the symmetric group. We give the result
for the  symmetric group of  degree $q^n$, where  $q$ is a power  of a
prime number $p$.

\begin{theorem}\cite[Theorem~4.8]{Cam}\label{th_cameron1}
Suppose $q$ is a power of the prime $p$.
\null

A maximal primitive subgroup $G$ of $\Sym(q^n)$ is one of the following:
\begin{enumerate}
\item\label{aff} affine, that is, $G = \AGL(d,p), p^d=q^n$, for some $d$;
\item\label{pnb} primitive  non-basic, that is, a wreath  product $G =
  \Sym(k)  \, \wr\, \Sym(r)$  in product  action, $k^r=q^n,  k\neq 2$,
  $r>1$.
\item\label{as} almost simple,  that is, $S \le G  \le \Aut(S)$, for a
  nonabelian simple group $S$.
\end{enumerate}

\end{theorem}  

Note that in our context $p = 2$.

It is convenient  to use a refinement of  the O'Nan-Scott theorem, due
to Cai  Heng Li~\cite{Li}, for the  special case when  $G$ contains an
abelian  regular subgroup  $T$;  in our  case,  this is  the group  of
translations.

\begin{theorem}\cite[Theorem~1.1]{Li}\label{th:Li}
  Let $G$ be a primitive group of degree $2^{d}$, with $d \ge 1$. Suppose
  $G$ contains a regular abelian subgroup $T$.  
  \null 
  
  Then $G$ is one of the following
  \begin{enumerate}
  \item\label{Liaff} affine, that is, $G \le \AGL(d,2)$;
  \item\label{Lipnb} 
    \begin{equation*}
      G = (S_{1} \times \dots \times S_{r}) . O . P,
    \end{equation*}
    with $2^{d}  = m^{r}$ for some  $m$ and $r  > 1$. Here $T  = T_{1}
    \times \dots \times T_{r}$, with $T_{i} < S_{i} \cong \Alt(m)$ for
    each $i$, $O \le \Out(S_{1}) \times \dots \times \Out(S_{r})$, and
    $P$ permutes transitively the $S_{i}$.
  \item\label{Lias} almost simple,  that is, $S \le G  \le \Aut(S)$, for a
    nonabelian simple group $S$.
  \end{enumerate}
\end{theorem}  

To prove the first statement of Theorem~\ref{uth} we need to deal with
the three possible cases of Theorem~\ref{th:Li}.

Case~\eqref{Liaff} is treated in Section~\ref{sec:affine}. An
important observation of Li \cite{Li} is in order here. If $V$ is a
vector space, with addition $+$, then the symmetric
group $\Sym(V)$  contains the affine group $\AGL(V) = T \GL(V)$, where
$T$ is the group of translations. But $\Sym(V)$ also contains
the conjugates of $\AGL(V)$, which are still affine groups on the
\emph{set} $V$, but possibly with respect to an operation $\circ$ different
from $+$. In particular the group $T$ of translations may be contained
in one of these conjugates, where it will be an abelian regular
subgroup. We have studied this situation in~\cite{Affine}, and we will be
exploiting these results in Section~\ref{sec:affine}.

Case~\eqref{Lipnb} will be dealt with in Section~\ref{sec:wreath}. 

In  the  almost  simple   case~\eqref{Lias},  the  intersection  of  a
one-point stabilizer  in $G$ with $S$  is a proper subgroup  of $S$ of
index  $2^{d}$,  since  the  nontrivial  normal subgroup  $S$  of  the
primitive group $G$ is transitive.  We can thus appeal (as Li does) to
a particular  case of a  result of Guralnick~\cite{Gur},  which states
that the only  nonabelian simple groups that have  a subgroup of index
of  the  form   $2^{d}$  are  either  the  alternating   groups  $S  =
\Alt(2^{d})$,  with  $d  >  2$,  or  the  groups  $\PSL(f,  q)$,
where $q$ is a prime-power, and $f$ is prime, 
$(q^{f}-1)/(q-1)= 2^{d}$.  We rule out the second possibility as
follows. Since $(q^{f}-1)/(q-1) = q^{f-1} + q^{f-2} + \dots +
q + 1 \equiv f \pmod{2}$, we have $f = 2$ here, and $q = 2^d -
1$. Well-known elementary arguments yield that $q$  and
$d$ are prime. However, $d = n_{t} m$ is not prime, as $n_{t} > 1$ by
assumption, and as noted earlier $m > 2$ by Cryptographic
Assumption~\eqref{item:01}.

Clearly $\Aut(\Alt(2^{d}))  = \Sym(2^{d})$ here, so $G$  is either the
alternating or the symmetric group. Since we have shown in
Section~\ref{sec:prelim} that $G \le \Alt(V)$, we obtain $G = \Alt(V)$.

To prove the second statement  of Theorem~\ref{uth}, we then appeal to
a standard argument:  if the nonabelian simple group  $G$ is generated
by a  subset $S$, then for  any fixed $r$  the set $S' =  \Set{s_1 s_2
  \dots s_r  : s_i  \in S}$  of $r$-fold products  of elements  of $S$
generates  a nontrivial  normal subgroup  of $G$,  and thus  $S'$ also
generates $G$. In our context $S$ is the set of the round functions
for all possible subkeys, and $r$
is the number of rounds, so that $S'$ is the set of the transformations
of the cryptosystem with independent subkeys.

\section{The affine case}
\label{sec:affine}

Suppose $G$  is contained in an  affine subgroup of  $\Sym(V)$. By the
theory  of~\cite{Affine},  there is  a  structure  of an  associative,
commutative, nilpotent ring  $(V, \circ, \cdot, 0)$ on  $V$, such that
$(V, \circ,  0)$ is a vector  space over the field  with two elements,
and ordinary addition on $V$ is expressed as
\begin{equation*}
  x + y = x \circ y \circ x y,
\end{equation*}
for  $x,  y  \in  V$.  Moreover,   $G$  acts  as  a  group  of  affine
transformations on $(V, \circ, 0)$.

As both $(V, \circ, 0)$ and $(V, +, 0)$ are elementary abelian, we have
\begin{equation*}
  0 = x + x = x \circ x \circ x x = 0 \circ x^{2} = x^{2}
\end{equation*}
for all $x \in V$. It follows
\begin{align*}
  x + y + x y 
  &= 
  (x \circ y \circ x y) \circ x y \circ (x \circ y \circ  x y) \cdot x
  y
  \\&=
  x \circ y \circ xy \circ xy \circ x^{2} y \circ x y^{2} \circ x^{2}
  y^{2}
  \\&=
  x \circ y.
\end{align*}
Here we have used the fact that $\cdot$ distributes over $\circ$.

Now $\rho \in G$ is linear  with respect to $\circ$, that is $(x \circ
y) \rho =  x\rho \circ y \rho$ for  all $x, y \in V$. Choose  $0 \ne y
\in U  = \Set{z \in V :  \text{$x z = 0$  for all $x \in  V$} }$. (The
latter set is different from $\{ 0 \}$, as the ring $(V, \circ, \cdot,
0)$ is nilpotent.) Then
\begin{equation}\label{eq:nearlin}
  (x + y) \rho = (x \circ y) \rho = x\rho \circ y \rho
  =
  x\rho + y \rho + x\rho \cdot y \rho.
\end{equation}

Now note that given $x \in V$, the set $x V = \Set{x z : z \in V}$ is
a subspace with respect to $\circ$, as $\cdot$ distributes over
$\circ$; and also a subspace with respect to $+$, as $x z_{1} + x z_{2} = x
z_{1} \circ x z_{2} \circ x^{2} z_{1} z_{2} = x
z_{1} \circ x z_{2}$.

It follows from~\ref{eq:nearlin} that for $0 \ne y \in U$ we have
\begin{equation*}
  \Set{(x + y)\rho + x \rho : x \in V}
  =
  y \rho + y \rho V.
\end{equation*}
The right hand side is a coset of a subspace of $V$ with respect to
$+$. Now $\lambda$ (and its inverse) are linear with respect to
$+$. Applying $\lambda^{-1}$ we obtain that
\begin{equation*}
  \Set{(x + y)\gamma + x \gamma : x \in V}
\end{equation*}
is also a coset of a subspace of $V$ with respect to
$+$. Choose an index $i$ so that the component $y_{i} \in V_{i}$ of
$y$ is nonzero. Then we have that the projection on $V_{i}$ of the
previous set
\begin{equation*}
  \Set{(x + y_{i})\gamma + x \gamma : x \in V_{i}}
\end{equation*}
is  a  coset of  a  subspace  of $V_{i}$  with  respect  to $+$.  This
contradicts Cryptographic Assumption~\eqref{item:r}.

\section{Wreath product in product action}
\label{sec:wreath}

Here we deal to the case when
\begin{equation*}
  G = (S_{1} \times \dots \times S_{r}) . O . P,
\end{equation*}
with $2^{d} = k^{r}$ for some $k$ and $r > 1$. Here $T = T_{1} \times
\dots \times T_{r}$, where $\Order{T_{i}} = k$ and $T_{i} < S_{i}
\cong \Alt(k)$ for each $i$, $O \le \Out(S_{1}) \times \dots \times
\Out(S_{r})$, and $P$ permutes transitively the $S_{i}$ by
conjugation. It follows that $S_{1} \times \dots \times S_{r} =
\Soc(G)$.

Note  that if  $k  =  2$ or  $4$,  so that  $S_{i}  \cong \Alt(2)$  or
$\Alt(4)$, the group $T$ of translations  is normal in $G$, so that $G
\le \AGL(V)$.  This contradicts  the non-linearity of  $\gamma$, which
follows  from Cryptographic  Assumption~\eqref{item:r}. Thus  we will
assume $k > 4$ in the rest of this section.

Note that  $G =  \Span{T, \rho}$, and  $T \le  \Soc(G)$, so that  $G /
\Soc(G)$ is cyclic, spanned by $\rho$. Since $P$ permutes transitively
the  $S_{i}$, it  follows that  $\rho$ permutes  \emph{cyclically} the
$S_{i}$  by  conjugation, that  is,  we  may  rename indices  so  that
$S_{i}^{\rho}  = \rho^{-1}  S_{i} \rho  = S_{i+1}$  for each  $i$ (and
indices are taken modulo $r$).

Since  each $T_{i}$  is  a group  of  translations, $W_{i}  = 0  T_{i}
\subseteq 0 S_{i}$ is a subspace of $V$, of order $k$. Since $0 S_{i}$
has also order $k$, $0 T_{i} = 0 S_{i}$. Clearly each element of $v \in
V$ can be written uniquely in the form $v = 0 t$, for $t \in T$. Thus
\begin{equation*}
 v = 0 t_{1} t_{2} \dots t_{r} = 0 t_{1} + 0 t_{2} + \dots + 0 t_{r}
\end{equation*}
for unique $t_{i} \in T_{i}$, and
\begin{equation*}
 V = W_{1} \oplus W_{2} \oplus \dots \oplus W_{r}.
\end{equation*}

For  each  $i$  we  have  also   $W_{i}  \rho  =  0  S_{i}  \rho  =  0
S_{i+1}^{\rho^{-1}} \rho = 0 \rho S_{i+1} = 0 S_{i+1} = W_{i+1}$,
as $0 \rho = 0$. Thus  $\rho$ permutes cyclically the $W_{i}$. Now let
$v \in V$, and write it as $v  = w_{1} + \dots + w_{r}$ where $w_i \in
W_i$. Let  $t_{i} \in  W_{i}$ be such  that $w_i  = 0 t_i$.  Since the
$t_{i}$ are translations, we  have $v = 0 t_{1} + 0  t_{2} + \dots + 0
t_{r} = 0  t_{1} t_{2} \dots t_{r}$.  We have $v \rho =  0 t_{1} t_{2}
\dots t_{r} \rho = 0 t_{1}^{\rho} t_{2}^{\rho} \dots  t_{r}^{\rho}$, as $0
\rho^{-1} = 0$. Since $t_{i}^{\rho}  \in S_{i}^{\rho} = S_{i+1}$, there are
$t_{i}'  \in T_{i}$  such that  $0  t_{i}^{\rho} =  0 t_{i}  \rho =  0
t_{i+1}'  \in  W_{i+1}$,  and  because  $S_{i}$  and  $S_{j}$  commute
elementwise, we have
\begin{align*}
  v \rho 
  &= 
  0 t_{1}^{\rho} t_{2}^{\rho} \dots  t_{r}^{\rho}
  =  0 t_{2}' t_{2}^{\rho} \dots t_{r}^{\rho} 
  =  0 t_{2}^{\rho} t_{2}'  \dots t_{r}^{\rho} 
  \\&=  0 t_{3}' t_{2}'  \dots t_{r}^{\rho} 
  =  0 t_{2}' t_{3}'  \dots t_{r}^{\rho} 
  = \dots
  \\&=  0 t_{2}' t_{3}'  \dots t_{1}' 
  = 0 t_{1}' + 0 t_{2}' + \dots + 0 t_{r}' 
  \\&= 0 t_{r} \rho + 0 t_{1} \rho + \dots + 0 t_{r-1} \rho
  \\&= w_{1} \rho + w_{2} \rho + \dots + w_{r} \rho.
\end{align*}

Now fix an index $i$, and  take $u \in W_{i}$.  We have from the above
\begin{equation*}
 v \rho = (w_{1} + w_{2} + \dots + w_{r})\rho 
 = w_{1} \rho + w_{2} \rho + +\cdots + w_{r} \rho,
\end{equation*}
where $w_i \rho \in W_{i+1}$, and also
\begin{equation*}
(v+u)\rho = w_{1} \rho + (w_{i} + u) \rho + \cdots +w_{r} \rho
\end{equation*}
with $(w_i+u)\rho \in W_{i+1}$. It follows
\begin{equation}\label{eq:difference}
(v+u)\rho + v\rho =  w_i\rho + (w_i + u)\rho \in W_{i+1}.
\end{equation}
Now  $\rho =  \gamma \lambda$,  where $\lambda$  is  linear.  Applying
$\lambda^{-1}$  to both  sides of~\eqref{eq:difference}  we get  $(v +
u)\gamma + v\gamma \in W_{i+1}\lambda^{-1}$. In other words, there are
subspaces $W_{i},  W_{i+1} \lambda^{-1}$ of $V$ of  the same dimension
such that when  the input difference to $\gamma$ is  in the first one,
then  the output  difference is  in second  one. By  the  arguments of
Section~\ref{sec:referee}  (with  $U  =   W_{i}$  and  $U'  =  W_{i+1}
\lambda^{-1}$), it follows  that $W_{i}$ is the direct  sum of some of
the  $V_{j}$,  for each  $i$.   Thus $W_{2}  =  W_{1}  \rho =  W_{1}
\lambda$ and $W_{3}  = W_{2} \lambda$, contradicting Cryptographic
Assumption~\eqref{item:affine}.


%
%

\providecommand{\bysame}{\leavevmode\hbox to3em{\hrulefill}\thinspace}
\providecommand{\MR}{\relax\ifhmode\unskip\space\fi MR }
\providecommand{\MRhref}[2]{%
  \href{http://www.ams.org/mathscinet-getitem?mr=#1}{#2}
}
\providecommand{\href}[2]{#2}

\end{document}